\pdfoutput=1
\RequirePackage{ifpdf}
\ifpdf % We are running pdfTeX in pdf mode
\documentclass[pdftex]{sigma}
\else
\documentclass{sigma}
\fi

\newtheorem{Theorem}{Theorem}[section]
\newtheorem{Corollary}[Theorem]{Corollary}
\newtheorem{Lemma}[Theorem]{Lemma}
\newtheorem{Proposition}[Theorem]{Proposition}
{\theoremstyle{definition}

\newtheorem{quest}[Theorem]{Question}
\newtheorem{Remark}[Theorem]{Remark}
\newtheorem{claim}{Claim}
}

\newcommand{\R}{\mathbb{R}}

\newcommand{\N}{\mathbb{N}}

\newcommand{\CC}{\mathcal{C}}
\newcommand{\DD}{\mathcal{D}}

\newcommand{\LL}{\mathcal{L}}
\newcommand{\PP}{\mathcal{P}}

\renewcommand{\SS}{\mathcal{S}}

\newcommand{\Cinf}{C^\infty}

\begin{document}

\allowdisplaybreaks

\renewcommand{\PaperNumber}{004}

\FirstPageHeading

\ShortArticleName{Embedding Theorems for the Dunkl Harmonic Oscillator on the Line}

\ArticleName{Embedding Theorems for the Dunkl Harmonic\\
 Oscillator on the Line}

\Author{Jes\'us A.~\'ALVAREZ L\'OPEZ~$^\dag$ and Manuel CALAZA~$^\ddag$}

\AuthorNameForHeading{J.A.~\'Alvarez L\'opez and M.~Calaza}

\Address{$^\dag$~Departamento de Xeometr\'{\i}a e Topolox\'{\i}a, Facultade de Matem\'aticas,\\
\hphantom{$^\dag$}~Universidade de Santiago de Compostela, 15782 Santiago de Compostela, Spain}
\EmailD{\href{mailto:jesus.alvarez@usc.es}{jesus.alvarez@usc.es}}

\Address{$^\ddag$~Laboratorio de Investigaci\'on 10, Servicio de Reumatolog\'ia, Instituto de Investigaci\'o³n\\
\hphantom{$^\ddag$}~Sanitaria, Hospital Cl\'inico Universitario, 15706 Santiago de Compostela, Spain}
\EmailD{\href{mailto:manuel.calaza@usc.es}{manuel.calaza@usc.es}}

\ArticleDates{Received September 09, 2013, in f\/inal form January 06, 2014; Published online January 10, 2014}

\Abstract{Embedding results of Sobolev type are proved for the Dunkl harmonic oscillator on the line.}

\Keywords{Dunkl harmonic oscillator; Sobolev embedding; generalized Hermite functions; Schwartz space}

\Classification{46E35; 47B25; 33C45}

\section{Introduction} \label{s:intro}

The subindex ev/odd is added to any space of functions on $\R$ to indicate its subspace of even/odd functions; in particular, $\Cinf=\Cinf_{\text{\rm ev}}\oplus\Cinf_{\text{\rm odd}}$ for $\Cinf:=\Cinf(\R)$. The Dunkl operator $T_\sigma$ ($\sigma>-1/2$) on $\Cinf$ is the perturbation of $\frac{{\rm d}}{{\rm d}x}$ def\/ined by $T_\sigma=\frac{{\rm d}}{{\rm d}x}$ on $\Cinf_{\text{\rm ev}}$ and $T_\sigma=\frac{{\rm d}}{{\rm d}x}+2\sigma\frac{1}{x}$ on $\Cinf_{\text{\rm odd}}$. The corresponding Dunkl harmonic oscillator is the perturbation $L_\sigma=-T_\sigma^2+s^2x^2$ of the harmonic oscillator $H=-\frac{{\rm d}^2}{{\rm d}x^2}+s^2x^2$ ($s>0$). The conjugation $E_\sigma=|x|^\sigma T_\sigma|x|^{-\sigma}$ on $|x|^\sigma\Cinf$ is equal to $\frac{{\rm d}}{{\rm d}x}-\sigma x^{-1}$ on $|x|^\sigma\Cinf_{\text{\rm ev}}$ and $\frac{{\rm d}}{{\rm d}x}+\sigma x^{-1}$ on $|x|^\sigma\Cinf_{\text{\rm odd}}$; note that $|x|^\sigma\Cinf_{\text{\rm ev/odd}}$ consists of even/odd functions, possibly not smooth or not even def\/ined at~$0$. Up to the product by a~constant, $E_\sigma$ was introduced by Yang~\cite{Yang1951}. In the form~$T_\sigma$, this operator was generalized to~$\R^n$ by Dunkl~\cite{Dunkl1988,Dunkl1989,Dunkl1991}, giving rise to what is now called Dunkl theory (see the survey~\cite{Rosler2003}); in particular, the Dunkl harmonic oscillator on~$\R^n$ was studied in \cite{Dunkl2002,NowakStempak2009b,NowakStempak2009a, Rosenblum1994}. See \cite{Plyushchay2000} for further generalizations on~$\R$. Sometimes the terms Yang--Dunkl operator and Yang--Dunkl harmonic oscillator are used in the case of~$\R$~\cite{Plyushchay2000}.

Let $p_k$ be the sequence of orthogonal polynomials for the measure $e^{-sx^2}|x|^{2\sigma}{\rm d}x$, taken with norm one and positive leading coef\/f\/icient. Up to normalization, these are the generalized Hermite polynomials~\cite[p.~380, Problem~25]{Szego1975}; see also \cite{Chihara1955,Chihara1978,DickinsonWarsi1963,DuttaChatterjeaMore1975,Rosenblum1994,Rosler1998}. The corresponding generalized Hermite functions are $\phi_k=p_ke^{-sx^2/2}$.

For each $m\in\N$, let $\SS^m$ be the Banach space of functions $\phi\in C^m(\R)$ with \mbox{$\sup_x|x^i\phi^{(j)}(x)|\!<\!\infty$} for $i+j\le m$; the corresponding Fr\'echet space $\SS=\bigcap_m\SS^m$ is the Schwartz space on $\R$. With domain $\SS$, $L_\sigma$ is essentially self-adjoint in $L^2(\R,|x|^{2\sigma}{\rm d}x)$, and the spectrum of its self-adjoint extension, $\LL_\sigma$, consists of the eigenvalues $(2k+1+2\sigma)s$ ($k\in\N$), with corresponding eigenfunctions $\phi_k$  \cite{Rosenblum1994}. For each real $m\ge0$, let $W_\sigma^m$ be the Hilbert space completion of $\SS$ with respect to the scalar product $\langle\phi,\psi\rangle_{W_\sigma^m}:=\langle(1+\LL_\sigma)^m\phi,\psi\rangle_\sigma$, where $\langle\ ,\ \rangle_\sigma$ denotes the scalar product of $L^2(\R,|x|^{2\sigma}{\rm d}x)$, obtaining a Fr\'echet space $W_\sigma^\infty=\bigcap_mW_\sigma^m$. We show the following embedding theorems; the second one is of Sobolev type.

\begin{Theorem}\label{t:SS^m' subset W_sigma^m}
  For each $m\in\N$, $\SS_{\text{\rm ev/odd}}^{M_{m,\text{\rm ev/odd}}}\subset W_{\sigma,\text{\rm ev/odd}}^m$ continuously, where
    \begin{gather*}
      M_{m,\text{\rm ev/odd}} =
        \begin{cases}
          \dfrac{3m+3}{2}+\dfrac{m+1}{4} \lceil\sigma\rceil(\lceil\sigma\rceil+3)+\lceil\sigma\rceil & \text{if $\sigma\ge0$ and $m$ is odd},\\
          2m+3 & \text{if $\sigma<0$ and $m$ is odd} ,
        \end{cases}
      \\
      M_{m,\text{\rm ev}} =
        \begin{cases}
          \dfrac{3m+2}{2}+\dfrac{m}{4} \lceil\sigma\rceil(\lceil\sigma\rceil+3)+\lceil\sigma\rceil & \text{if $\sigma\ge0$ and $m$ is even},\\
          2m+2 & \text{if $\sigma<0$ and $m$ is even} ,
        \end{cases}
      \\
      M_{m,\text{\rm odd}} =
        \begin{cases}
          \dfrac{3m+4}{2}+\dfrac{m+2}{4} \lceil\sigma\rceil(\lceil\sigma\rceil+3)+\lceil\sigma\rceil & \text{if $\sigma\ge0$ and $m$ is even},\\
          2m+4 & \text{if $\sigma<0$ and $m$ is even}.
        \end{cases}
    \end{gather*}
\end{Theorem}

\begin{Theorem}\label{t:Sobolev}
  For all $m\in\N$ and $m_\sigma=m+1+\frac{1}{2}\lceil\sigma\rceil(\lceil\sigma\rceil+1)$, $W_{\sigma,\text{\rm ev/odd}}^{m'}\subset\SS_{\text{\rm ev/odd}}^m$ continuously if $m'>N_{m,\text{\rm ev/odd}}$, where $N_{m,\text{\rm ev}}=2$ and $N_{m,\text{\rm odd}}=5$ if $m_\sigma=1$, $N_{m,\text{\rm ev}}=6$ and $N_{m,\text{\rm odd}}=5$ if $m_\sigma=2$, $N_{m,\text{\rm ev}}=6$ and $N_{m,\text{\rm odd}}=7$ if $m_\sigma=3$, and $N_{m,\text{\rm ev/odd}}=m_\sigma+3$ for $m_\sigma\ge4$.
\end{Theorem}

\begin{Corollary}\label{c:SS=W_sigma^infty}
  $\SS=W_\sigma^\infty$ as Fr\'echet spaces.
\end{Corollary}

In other words, Corollary~\ref{c:SS=W_sigma^infty} states that an element $\phi\in L^2(\R,|x|^{2\sigma}{\rm d}x)$ is in $\SS$ if and only if the ``Fourier coef\/f\/icients'' $\langle\phi,\phi_k\rangle_\sigma$ are rapidly decreasing on $k$. This also means that $\SS=\bigcap_m\DD(\LL_\sigma^m)$ ($\SS$~is the smooth core of $\LL_\sigma$ with the terminology of \cite{BruningLesch1992}) because the sequence of eigenvalues of~$\LL_\sigma$ is in~$O(k)$ as~$k\to\infty$.

We introduce a version $\SS_\sigma^m$ of every $\SS^m$, whose def\/inition involves $T_\sigma$ instead of $\frac{{\rm d}}{{\rm d}x}$. They satisfy much simpler embeddings: $S_\sigma^{\lceil m\rceil+1}\subset W_\sigma^m$, and $W_\sigma^{m'}\subset\SS_\sigma^m$ if $m'-m>1$. Even though $\SS=\bigcap_m\SS_\sigma^m$, the inclusion relations between the spaces $\SS_\sigma^m$ and $\SS^{m'}$ are complicated, giving rise to the complexity of Theorems~\ref{t:SS^m' subset W_sigma^m} and~\ref{t:Sobolev}.

Other Sobolev type embedding theorems, for dif\/ferent operators and with dif\/ferent techniques, were recently proved in \cite{Watanabe1997,Watanabe2002,Watanabe2004}.

Next, we consider other perturbations of $H$ on $\R_+$. Let $\SS_{\text{\rm ev},U}$ denote the space of restrictions of even Schwartz functions to some open $U\subset\R_+$, and set $\phi_{k,U}=\phi_k|_U$.

\begin{Theorem}\label{t:P}
  Let $P=H-2f_1 \frac{{\rm d}}{{\rm d}x}+f_2$, where $f_1\in C^1(U)$ and $f_2\in C(U)$ for some open $U\subset\R_+$ of full Lebesgue measure. Assume that $f_2=\sigma(\sigma-1)x^{-2}-f_1^2-f_1'$ for some $\sigma>-1/2$. Let $h=x^\sigma e^{-F_1}$, where $F_1\in C^2(U)$ is a primitive of $f_1$. Then the following properties hold:
    \begin{itemize}\itemsep=0pt

      \item[$(i)$] $P$, with domain $h \SS_{\text{\rm ev},U}$, is essentially self-adjoint in $L^2(U,e^{2F_1}{\rm d}x)$;

      \item[$(ii)$] the spectrum of its self-adjoint extension, $\PP$, consists of the eigenvalues $(4k+1+2\sigma)s$ {\rm(}$k\in\N${\rm)} with multiplicity one and normalized eigenfunctions $\sqrt{2} h\phi_{2k,U}$; and

       \item[$(iii)$] the smooth core of $\PP$ is $h\SS_{\text{\rm ev},U}$.

    \end{itemize}
\end{Theorem}

This theorem follows by showing that the stated condition on $f_1$ and $f_2$ characterizes the cases where $P$ can be obtained by the following process: f\/irst, restricting $L_\sigma$ to even functions, then restricting to~$U$, and f\/inally conjugating by $h$. The term of $P$ with $\frac{{\rm d}}{{\rm d}x}$ can be removed by conjugation with the product of a positive function, obtaining the operator $H+\sigma(\sigma-1)x^{-2}$; in this way, we get all operators of the form $H+cx^{-2}$ with $c>-1/4$.

The conditions of Theorem~\ref{t:P} are satisf\/ied by $P=H-2c_1x^{-1}\frac{{\rm d}}{{\rm d}x}+c_2x^{-2}$ ($c_1,c_2\in\R$) on~$\R_+$ if and only if there is some $a\in\R$ such that $a^2+(2c_1-1)a-c_2=0$ and $a+c_1>-1/2$; in this case, $h=x^a$ and $e^{2F_1}=x^{2c_1}$. For some $c_1,c_2\in\R$, there are two values of~$a$ satisfying these conditions, obtaining two dif\/ferent self-adjoint operators def\/ined by $P$ in dif\/ferent Hilbert spaces. For instance, $L_\sigma$ may def\/ine a self-adjoint operator when~$\sigma\le-1/2$.

This example is applied in \cite{AlvCalaza:Witten} to prove a new type of Morse inequalities on strata of compact stratif\/ications \cite{Mather1970,Thom1969,Verona1984} with adapted metrics~\cite{BrasseletHectorSaralegi1992, Nagase1983,Nagase1986}, where the Witten's perturbation \cite{Witten1982} is used for the minimal/maximal ideal boundary conditions of de~Rham complex \cite{BruningLesch1992, Cheeger1980,Cheeger1983}. The version of Morse functions used in \cite{AlvCalaza:Witten} is dif\/ferent from the version of Goresky--MacPherson \cite{GoreskyMacPherson1983:Morse}. More precisely, in the local conic model of a stratif\/ication around each critical point, by induction on the depth of the stratif\/ication, it is assumed that the Laplacian of the minimal/maximal boundary condition of the de~Rham complex of the link (section of the cone) has a nice spectral decomposition. Using this, the Witten's perturbation of the de~Rham complex of the cone splits into an inf\/inite direct sum of elliptic complexes of two simple types, with length one and two, which represent the radial direction of Witten's perturbed complex. It turns out that the Witten's perturbed Laplacian of these simple complexes is described by the above operator $P$, and the two possible choices of the constant~$a$ give rise to the minimal/maximal ideal boundary conditions. In this way, Theorem~\ref{t:P} becomes a key ingredient of~\cite{AlvCalaza:Witten}.

\section{Preliminaries}\label{s:preliminaries}

\subsection{Dunkl operator on the line}\label{ss:D sigma}

For any $\phi\in\Cinf:=\Cinf(\R)$, there exists some $\psi\in\Cinf$ so that $\phi(x)-\phi(0)=x\psi(x)$; moreover
  \begin{gather}\label{psi (m)(x)}
    \psi^{(m)}(x)=\int_0^1t^m\phi^{(m+1)}(tx){\rm d}t
  \end{gather}
for all\footnote{We adopt the convention $0\in\N$.} $m\in\N$ (see e.g.~\cite[Theorem~1.1.9]{Hormander1990-I}). Let us use the notation $\psi=x^{-1}\phi$. The Dunkl operator on $T_\sigma$ ($\sigma\in\R$) on $\Cinf$ is the perturbation of $\frac{{\rm d}}{{\rm d}x}$ def\/ined by
  \[
    (T_\sigma\phi)(x)=\phi'(x)+2\sigma \frac{\phi(x)-\phi(-x)}{x} .
  \]

Consider matrix expressions of operators on $\Cinf$ with respect to the decomposition $\Cinf=\Cinf_{\text{\rm ev}}\oplus\Cinf_{\text{\rm odd}}$, as direct sum of subspaces of even and odd functions. For each function $h$, the notation $h$ is also used for the operator of multiplication by $h$. Then
  \begin{gather*}
    \frac{{\rm d}}{{\rm d}x}=
      \begin{pmatrix}
        0 &  \dfrac{{\rm d}}{{\rm d}x} \\
        \dfrac{{\rm d}}{{\rm d}x} & 0
      \end{pmatrix} ,\qquad
    x=
      \begin{pmatrix}
        0 &  x \\
        x & 0
      \end{pmatrix} ,\\
    T_\sigma=
      \begin{pmatrix}
        0 &  \dfrac{{\rm d}}{{\rm d}x}+2\sigma x^{-1} \\
        \dfrac{{\rm d}}{{\rm d}x} & 0
      \end{pmatrix}
    = \frac{{\rm d}}{{\rm d}x}+2\sigma
      \begin{pmatrix}
        0 &  x^{-1} \\
        0 & 0
      \end{pmatrix}
  \end{gather*}
on $\Cinf$. With
$
  \Sigma=
    \begin{pmatrix}
      \sigma & 0 \\
      0 & -\sigma
    \end{pmatrix}
$,
we get
  \begin{gather}\label{[D sigma,x]=1+2Sigma, ...}
    [T_\sigma,x]=1+2\Sigma ,\qquad
    T_\sigma\Sigma+\Sigma T_\sigma=x \Sigma+\Sigma x=0 .
  \end{gather}

\subsection{Dunkl harmonic oscillator on the line}\label{ss:L}

On $\Cinf$, the harmonic oscillator, and the annihilation and creation operators are $H=-\frac{{\rm d}^2}{{\rm d}x^2}+s^2x^2$, $A=sx+\frac{{\rm d}}{{\rm d}x}$ and $A'=sx-\frac{{\rm d}}{{\rm d}x}$ ($s>0$). Their perturbations $L=-T_\sigma^2+s^2x^2$, $B=sx+T_\sigma$ and $B'=sx-T_\sigma$ are called Dunkl harmonic oscillator, and Dunkl annihilation and creation operators. By~\eqref{[D sigma,x]=1+2Sigma, ...},
  \begin{gather}
    L=BB'-(1+2\Sigma)s=B'B+(1+2\Sigma)s=\frac{1}{2}(BB'+B'B) ,\label{L}\\
    [L,B]=-2sB ,\qquad[L,B']=2sB' ,\label{[L,B]}\\
    [B,B']=2s(1+2\Sigma) ,\label{[B,B']}\\
    [L,\Sigma]=B\Sigma+\Sigma B=B'\Sigma+\Sigma B'=0 .
    \label{[L,Sigma]=B' Sigma+Sigma B'=0}
  \end{gather}

For each $m\in\N$, let $\SS^m$ be the space of functions $\phi\in\Cinf$ such that
  \[
    \|\phi\|_{\SS^m}=\sum_{i+j\le m}\sup_x|x^i\phi^{(j)}(x)|<\infty .
  \]
This expression def\/ines a norm $\|\; \|_{\SS^m}$ on $\SS^m$, which becomes a Banach space. We have \mbox{$\SS^{m+1}\subset\SS^m$} continuously\footnote{For topological vector spaces $X$ and $Y$, it is said that $X\subset Y$ continuously if $X$ is a linear subspace of $Y$ and the inclusion map $X\hookrightarrow Y$ is continuous.}, and $\SS=\bigcap_m\SS^m$, with the induced Fr\'echet topology, is the Schwartz space on~$\R$. Note that $\|\phi'\|_{\SS^m}\le\|\phi\|_{\SS^{m+1}}$ for all~$m$.

We can restrict the above decomposition of $\Cinf$ to every $\SS^m$ and $\SS$,
obtaining $\SS^m=\SS^m_{\text{\rm ev}}\oplus\SS^m_{\text{\rm odd}}$ and $\SS=\SS_{\text{\rm ev}}\oplus\SS_{\text{\rm odd}}$. The matrix expressions of operators on $\SS$ are taken with respect to this decomposition. For $\phi\in \Cinf_{\text{\rm ev}}$, $\psi=x^{-1}\phi$ and $i,j\in\N$, we get from~\eqref{psi (m)(x)} that
  \[
    \big|x^i\psi^{(j)}(x)\big|\le\int_0^1t^{j-i}\big|(tx)^i\phi^{(j+1)}(tx)\big|{\rm d}t\le\sup_{y\in\R}\big|y^i\phi^{(j+1)}(y)\big|
  \]
for all $x\in\R$. So $\|\psi\|_{\SS^m}\le\|\phi\|_{\SS^{m+1}}$ for all $m\in\N$, obtaining that $\SS_{\text{\rm odd}}=x \SS_{\text{\rm ev}}$ and $x^{-1}:\Cinf_{\text{\rm odd}}\to \Cinf_{\text{\rm ev}}$ restricts to a continuous operator $x^{-1}:\SS_{\text{\rm odd}}\to\SS_{\text{\rm ev}}$. Hence $x:\SS_{\text{\rm ev}}\to\SS_{\text{\rm odd}}$ is an isomorphism of Fr\'echet spaces, and $T_\sigma$, $B$, $B'$ and $L$ def\/ine continuous operators on $\SS$.

Let $\langle\ ,\ \rangle_\sigma$ and $\|\ \|_\sigma$ be the scalar product and norm of $L^2(\R,|x|^{2\sigma}{\rm d}x)$. Suppose from now on that $\sigma>-1/2$, obtaining that $\SS$ is dense in $L^2(\R,|x|^{2\sigma}{\rm d}x)$. The following properties hold considering these operators in $L^2(\R,|x|^{2\sigma}{\rm d}x)$ with domain $\SS$: $-T_\sigma$ is adjoint of $T_\sigma$, $B'$ is adjoint of $B$, and $L$ is essentially self-adjoint. Let $\LL$, or $\LL_\sigma$, denote the self-adjoint extension of~$L$ (with domain~$\SS$). Its spectrum consists of the eigenvalues $(2k+1+2\sigma)s$ ($k\in\N$). The corresponding normalized eigenfunctions $\phi_k$ are inductively def\/ined by
  \begin{gather}
    \phi_0 =s^{(2\sigma+1)/4}\Gamma(\sigma+1/2)^{-1/2}e^{-sx^2/2} ,\label{phi 0}\\
    \phi_k =
      \begin{cases}
        (2ks)^{-1/2}B'\phi_{k-1} & \text{if $k$ is even},\\
        (2(k+2\sigma)s)^{-1/2}B'\phi_{k-1} & \text{if $k$ is odd},
      \end{cases}
    \qquad k\ge1 .
    \label{phi k}
  \end{gather}
Furthermore
    \begin{gather}
      B\phi_0 =0 ,\label{B phi 0}\\
      B\phi_k =
        \begin{cases}
          (2ks)^{1/2}\phi_{k-1} & \text{if $k$ is even},\\
          (2(k+2\sigma)s)^{1/2}\phi_{k-1} & \text{if $k$ is odd},
        \end{cases}
      \qquad k\ge1 .
      \label{B phi k}
    \end{gather}
These properties follow from \eqref{L}--\eqref{[L,Sigma]=B' Sigma+Sigma B'=0}, like in the case of $H$.

\subsection{Generalized Hermite polynomials}\label{ss:p k}

By~\eqref{phi 0},~\eqref{phi k} and the def\/inition of $B'$, we get $\phi_k=p_ke^{-sx^2/2}$, where $p_k$ is the sequence of polynomials inductively given by $p_0=s^{(2\sigma+1)/4}\Gamma(\sigma+1/2)^{-1/2}$ and
  \begin{gather}\label{p k}
    p_k=
      \begin{cases}
        (2ks)^{-1/2}(2sxp_{k-1}-T_\sigma p_{k-1}) & \text{if $k$ is even},\\
        (2(k+2\sigma)s)^{-1/2}(2sxp_{k-1}-T_\sigma p_{k-1}) & \text{if $k$ is odd},
      \end{cases}
    \qquad k\ge1 .
  \end{gather}
Up to normalization, $p_k$ and $\phi_k$ are the generalized Hermite polynomials and functions \cite[p.~380, Problem~25]{Szego1975}. Each $p_k$ is of degree~$k$, even/odd if~$k$ is even/odd, and with positive leading coef\/f\/icient. Moreover $T_\sigma p_0=0$ and
  \begin{gather}\label{D sigma p k}
    T_\sigma p_k=
          \begin{cases}
            (2ks)^{1/2}p_{k-1} & \text{if $k$ is even},\\
            (2(k+2\sigma)s)^{1/2}p_{k-1} & \text{if $k$ is odd},
          \end{cases}
        \qquad k\ge1 .
  \end{gather}
From~\eqref{p k} and~\eqref{D sigma p k}, we obtain the recursion formula
   \begin{gather}\label{recurrence}
     p_k=
       \begin{cases}
         k^{-1/2}\big((2s)^{1/2}xp_{k-1}-(k-1+2\sigma)^{1/2}p_{k-2}\big) & \text{if $k$ is even},\\
         (k+2\sigma)^{-1/2}\big((2s)^{1/2}xp_{k-1}-(k-1)^{1/2}p_{k-2}\big) & \text{if $k$ is odd} .
       \end{cases}
   \end{gather}
By~\eqref{recurrence} and induction on $k$, we easily get the following when $k$ is odd\footnote{As a convention, the product of an empty set of factors is~$1$. Thus $(k-1)(k-3)\cdots(\ell+2)=1$ for $\ell=k-1$ in~\eqref{x -1 p k}.}:
   \begin{gather}\label{x -1 p k}
     x^{-1}p_k=\sum_{\ell\in\{0,2,\dots,k-1\}}(-1)^{\frac{k-\ell-1}{2}}
     \sqrt{\frac{(k-1)(k-3)\cdots(\ell+2)2s}{(k+2\sigma)(k-2+2\sigma)\cdots(\ell+1+2\sigma)}} p_\ell .
  \end{gather}

The following theorem contains a simplif\/ied version of the asymptotic estimates satisf\/ied by~$\phi_k$ and $\xi_k=|x|^\sigma\phi_k$. They can be obtained by expressing~$p_n$ in terms of the Laguerre  polynomials~\cite{Rosler1998,Rosler2003}, whose asymptotic estimates are studied in \cite{AskeyWainger1965,Erdelyi1960,Muckenhoupt1970b, Muckenhoupt1970a-II}. The method of Bonan--Clark~\cite{BonanClark1990} can be also applied~\cite{AlvCalaza:generalizedHermite}.

  \begin{Theorem}\label{t:upper estimates of xi k}
  There exist $C,C',C''>0$, depending on $\sigma$ and $s$, such that:
    \begin{itemize}\itemsep=0pt

      \item[$(i)$] if $k$ is odd or $\sigma\ge0$, then $\xi_k^2(x)\le C'k^{-1/6}$ for all $x\in\R$;

      \item[$(ii)$] if $k$ is even and positive, and $\sigma<0$, then $\xi_k^2(x)\le C''k^{-1/6}$ for $|x|\ge 1$; and

      \item[$(iii)$] if $\sigma<0$, then $\phi_k^2(x)\le C''$ for all $k$ and $|x|\le1$.

    \end{itemize}
\end{Theorem}

\section{Perturbed Schwartz space}\label{s:SS sigma}

We introduce a perturbed version $\SS_\sigma^m$ of each $\SS^m$ that will be appropriate to show our embedding results. Since~$\SS_\sigma^m$ must contain the functions $\phi_k$, Theorem~\ref{t:upper estimates of xi k} indicates that dif\/ferent def\/initions must be given for $\sigma\ge0$ and $\sigma<0$.

When $\sigma\ge0$, for any $\phi\in\Cinf$ and $m\in\N$, let
  \begin{gather}\label{SS_sigma^m, sigma ge0}
    \|\phi\|_{\SS_\sigma^m}=\sum_{i+j\le m}\sup_x|x|^\sigma |x^iT_\sigma^j\phi(x)| .
  \end{gather}
This def\/ines a norm $\|\ \|_{\SS_\sigma^m}$ on the linear space of
functions $\phi\in\Cinf$ with $\|\phi\|_{\SS_\sigma^m}<\infty$,
and let~$\SS_\sigma^m$ denote the corresponding Banach space
completion. There are direct sum decompositions into
subspaces of even and odd functions, $\SS_\sigma^m=\SS_{\sigma,\text{\rm ev}}^m\oplus\SS_{\sigma,\text{\rm odd}}^m$.

When $\sigma<0$, the even and odd functions are considered separately: let
  \begin{gather}
    \|\phi\|_{\SS_\sigma^m}
    =\sum_{i+j\le m,\ i+j\ \text{even}}\left(\sup_{|x|\le1}|x^i(T_\sigma^j\phi)(x)|
    +\sup_{|x|\ge1}|x|^\sigma |x^i(T_\sigma^j\phi)(x)|\right)\nonumber\\
\hphantom{\|\phi\|_{\SS_\sigma^m}=}{}+\sum_{i+j\le m,\ i+j\ \text{odd}}\sup_{x\ne0}|x|^\sigma |x^i(T_\sigma^j\phi)(x)|\label{SS_sigma,even^m, sigma<0}
   \end{gather}
for $\phi\in\Cinf_{\text{\rm ev}}$, and
  \begin{gather}
    \|\phi\|_{\SS_\sigma^m}
    =\sum_{i+j\le m,\ i+j\ \text{even}}\sup_{x\ne0}|x|^\sigma |x^i(T_\sigma^j\phi)(x)|\nonumber\\
\hphantom{\|\phi\|_{\SS_\sigma^m}=}{}+\sum_{i+j\le m,\ i+j\ \text{odd}}\left(\sup_{|x|\le1}|x^i(T_\sigma^j\phi)(x)|
    +\sup_{|x|\ge1}|x|^\sigma |x^i(T_\sigma^j\phi)(x)|\right)\label{SS_sigma,odd^m, sigma<0}
  \end{gather}
for $\phi\in\Cinf_{\text{\rm odd}}$. These expressions def\/ine a norm
$\|\ \|_{\SS_\sigma^m}$ on the linear spaces of
functions $\phi\in\Cinf_{\text{\rm ev/odd}}$ with
$\|\phi\|_{\SS_\sigma^m}<\infty$. The corresponding
Banach space completions will be denoted by $\SS_{\sigma,\text{\rm
ev/odd}}^m$, and let $\SS_\sigma^m=\SS_{\sigma,\text{\rm ev}}^m\oplus\SS_{\sigma,\text{\rm odd}}^m$.

In any case, there are continuous inclusions
$\SS_\sigma^{m+1}\subset\SS_\sigma^m$, and a perturbed Schwartz space is def\/ined as $\SS_\sigma=\bigcap_m\SS_\sigma^m$, with the corresponding Fr\'echet topology, which decomposes as direct sum of the subspaces of even and odd functions, $\SS_\sigma=\SS_{\sigma,\text{\rm ev}}\oplus\SS_{\sigma,\text{\rm odd}}$; in particular, $\SS_0=\SS$.  It easily follows
that $\SS_\sigma$ consists of functions that are $\Cinf$ on $\R\setminus\{0\}$ but \textit{a priori} possibly not even def\/ined at zero, and $\SS_\sigma^m\cap\Cinf$ is dense in $\SS_\sigma^m$ for all $m$; thus $\SS_\sigma\cap\Cinf$ is dense in $\SS_\sigma$.

Obviously, $\Sigma$ def\/ines a bounded operator on each $\SS_\sigma^m$. It is also easy to see that $T_\sigma$ def\/ines a~bounded operator $\SS_\sigma^{m+1}\to\SS_\sigma^m$ for any $m$; notice that, when $\sigma<0$, the role played by the parity of $i+j$ f\/its well to prove this property. Similarly, $x$ def\/ines a bounded operator $\SS_\sigma^{m+1}\to\SS_\sigma^m$ for any $m$ because, by~\eqref{[D sigma,x]=1+2Sigma, ...},
  \[
    [T_\sigma^j,x]=
      \begin{cases}
        jT_\sigma^{j-1} & \text{if $j$ is even},\\
        (j+2\Sigma)T_\sigma^{j-1} & \text{if $j$ is odd}.
      \end{cases}
  \]
So $B$ and $B'$ def\/ine bounded operators $\SS_\sigma^{m+1}\to\SS_\sigma^m$, and $L$ a bounded operator $\SS_\sigma^{m+2}\to\SS_\sigma^m$. Thus $T_\sigma$, $x$, $\Sigma$, $B$, $B'$ and $L$ def\/ine continuous operators on $\SS_\sigma$.

In order to prove Theorems~\ref{t:SS^m' subset W_sigma^m} and~\ref{t:Sobolev}, we introduce an intermediate
weakly perturbed Schwartz space $\SS_{w,\sigma}$. Like $\SS_\sigma$, it is a Fr\'echet space of the form $\SS_{w,\sigma}=\bigcap_m\SS_{w,\sigma}^m$, where each $\SS_{w,\sigma}^m$ is the Banach space def\/ined like $\SS_\sigma^m$ by using  $\frac{{\rm d}}{{\rm d}x}$ instead of $T_\sigma$ in the right hand sides of~\mbox{\eqref{SS_sigma^m, sigma ge0}--\eqref{SS_sigma,odd^m, sigma<0}}; in particular, $\SS_{w,\sigma}^0=\SS_\sigma^0$ as Banach spaces. Let $\|\ \|_{\SS_{w,\sigma}^m}$ denote the norm of $\SS_{w,\sigma}^m$. As before, $\SS_{w,\sigma}$ consists of functions which are $\Cinf$ on $\R\setminus\{0\}$ but \textit{a priori} possibly not even def\/ined at zero, $\SS_{w,\sigma}\cap\Cinf$ is dense in $\SS_{w,\sigma}$, there is a canonical decomposition $\SS_{w,\sigma}=\SS_{w,\sigma,\text{\rm ev}}\oplus\SS_{w,\sigma,\text{\rm odd}}$, and $\frac{{\rm d}}{{\rm d}x}$ and $x$ def\/ine bounded operators on $\SS_{w,\sigma}^{m+1}\to\SS_{w,\sigma}^m$. Thus  $\frac{{\rm d}}{{\rm d}x}$ and $x$ def\/ine continuous operators on $\SS_{w,\sigma}$.

\begin{Lemma}\label{l:SS m+lceil sigma/2 rceil subset SS w,sigma m}
  If $\sigma\ge0$, then $\SS^{m+\lceil\sigma\rceil}\subset\SS_{w,\sigma}^m$ continuously.
\end{Lemma}

\begin{proof}
  Let $\phi\in\SS$. For all $i$ and $j$, we have $|x|^\sigma|x^i\phi^{(j)}(x)|\le|x^{i+\lceil\sigma\rceil}\phi^{(j)}(x)|$ for $|x|\ge1$, and $|x|^\sigma|x^i\phi^{(j)}(x)|\le|x^i\phi^{(j)}(x)|$ for $|x|\le1$. So $\|\phi\|_{\SS_{w,\sigma}^m}\le\|\phi\|_{\SS^{m+\lceil\sigma\rceil}}$.
\end{proof}

\begin{Lemma}\label{l:SS_w,sigma^m' subset SS^m}
  If $\sigma\ge0$, $\SS_{w,\sigma}^{m_\sigma}\subset\SS^m$ continuously, where $m_\sigma=m+1+\frac{1}{2}\lceil\sigma\rceil(\lceil\sigma\rceil+1)$.
\end{Lemma}

\begin{proof}
  Let $\phi\in\SS_{w,\sigma}$. For all $i$ and $j$,
    \begin{gather}\label{|x| ge1}
      \big|x^i\phi^{(j)}(x)\big|\le |x|^\sigma\big|x^i\phi^{(j)}(x)\big|
    \end{gather}
  for $|x|\ge1$. It remains to prove an inequality of this type for $|x|\le1$, which will be a consequence of the following assertion.

  \begin{claim}\label{cl:phi(x)}
      For each $n\in\N$, there are f\/inite families of real numbers, $c^n_{a,b}$, $d^n_{k,\ell}$ and $e^n_{u,v}$, where the indices $a$, $b$, $k$, $\ell$, $u$ and $v$ run in f\/inite subsets of $\N$ with $b,\ell,v\le M_n=1+\frac{n(n+1)}{2}$ and $k\ge n$, such that, for all $\phi\in\Cinf$,
        \[
          \phi(x)=\sum_{a,b}c^n_{a,b}x^a\phi^{(b)}(1)+\sum_{k,\ell}d^n_{k,\ell}x^k\phi^{(\ell)}(x)
          +\sum_{u,v}e^n_{u,v}x^u\int_x^1t^n\phi^{(v)}(t){\rm d} t .
        \]
    \end{claim}

  Assuming that Claim~\ref{cl:phi(x)} is true, the proof can be completed as follows. Let $\phi\in\SS_{w,\sigma}$ and set $n=\lceil\sigma\rceil$. For $|x|\le1$, according to Claim~\ref{cl:phi(x)}, $|\phi(x)|$ is bounded by
    \begin{gather*}
      \sum_{a,b}|c^n_{a,b}|\big|\phi^{(b)}(1)\big|+\sum_{k,\ell}|d^n_{k,\ell}|\big|x^k\phi^{(\ell)}(x)\big|
       +\sum_{u,v}|e^n_{u,v}| 2\max_{|t|\le1}\big|t^n\phi^{(v)}(t)\big|\\
     \qquad{}  \le\sum_{i,j}|c^n_{a,b}|\big|\phi^{(b)}(1)\big|+\sum_{k,\ell}|d^n_{k,\ell}| |x|^\sigma\big|\phi^{(\ell)}(x)\big|
       +\sum_{u,v}|e^n_{u,v}| 2\max_{|t|\le1}|t|^\sigma\big|\phi^{(v)}(t)\big| .
    \end{gather*}
  Let $m,i,j\in\N$ with $i+j\le m$. By applying the above inequality to the function $x^i\phi^{(j)}$, and expressing each derivative $(x^i\phi^{(j)})^{(r)}$ as a linear combination of functions of the form $x^p\phi^{(q)}$ with $p+q\le i+j+r$, it follows that there is some $C\ge1$, depending only on $\sigma$ and $m$, such that
    \begin{gather}\label{|x| le1}
      \big|x^i\phi^{(j)}(x)\big|\le C \|\phi\|_{\SS_{w,\sigma}^{i+j+M_n}}
    \end{gather}
  for $|x|\le1$. By~\eqref{|x| ge1} and~\eqref{|x| le1}, $\|\phi\|_{\SS^m}\le C \|\phi\|_{\SS_{w,\sigma}^{m_\sigma}}$ because $m_\sigma=m+M_n$.

  Now, let us prove Claim~\ref{cl:phi(x)}. By induction on $n$ and using integration by parts, it is easy to prove that
    \begin{gather}\label{int x 1t n f (n+1)(t) dt}
      \int_x^1t^n\phi^{(n+1)}(t){\rm d}t=\sum_{r=0}^n(-1)^{n-r} \frac{n!}{r!} \bigl(\phi^{(r)}(1)-x^r\phi^{(r)}(x)\bigr) .
      \end{gather}
  This shows directly Claim~\ref{cl:phi(x)} for $n\in\{0,1\}$. Proceeding by induction, let $n\ge2$ and assume that Claim~\ref{cl:phi(x)} holds for $n-1$. By~\eqref{int x 1t n f (n+1)(t) dt}, it is enough to f\/ind appropriate expressions of $x^r\phi^{(r)}(x)$ for $0<r<n$. For that purpose, apply Claim~\ref{cl:phi(x)} for $n-1$ to each function $\phi^{(r)}$, and multiply the resulting equality by $x^r$ to get
    \begin{gather*}
    \begin{split}
    &  x^r\phi^{(r)}(x)=\sum_{a,b}c^{n-1}_{a,b}x^{r+a}\phi^{(r+b)}(1)+\sum_{k,\ell}d^{n-1}_{k,\ell}x^{r+k}\phi^{(r+\ell)}(x)\\
& \hphantom{x^r\phi^{(r)}(x)=}{}+\sum_{u,v}e^{n-1}_{u,v}x^{r+u}\int_x^1t^{n-1}\phi^{(r+v)}(t){\rm d}t ,
\end{split}
    \end{gather*}
  where $a$, $b$, $k$, $\ell$, $u$ and $v$ run in f\/inite subsets of $\N$ with $b,\ell,v\le M_{n-1}$ and $k\ge n-1$; thus $r+k\ge n$ and $r+b,r+\ell,r+v\le n-1+M_{n-1}=M_n-1$. Therefore it only remains to rise the exponent of $t$ by a unit in the integrals of the last sum. Once more, integration by parts makes the job:
    \begin{gather*}
      \int_x^1t^n\phi^{(r+v+1)}(t){\rm d}t=\phi^{(r+v)}(1)-x^n\phi^{(r+v)}(x)-n\int_x^1t^{n-1}\phi^{(r+v)}{\rm d}t .\tag*{\qed}
    \end{gather*}
\renewcommand{\qed}{}
\end{proof}

\begin{Lemma}\label{l:SS m+1 subset SS w,sigma m}
  If $\sigma<0$, then $\SS^{m+1}\subset\SS_{w,\sigma}^m$ continuously.
\end{Lemma}

\begin{proof}
  This is proved by induction on $m$. For $\phi\in C^\infty_{\text{\rm ev}}$ and $|x|\ge1$, we have $|x|^\sigma |\phi(x)|\le|\phi(x)|$, obtaining $\|\ \|_{\SS_{w,\sigma}^0}\le\|\ \|_{\SS^0}$ on $\Cinf_{\text{\rm ev}}$. On the other hand, for $\phi\in\Cinf_{\text{\rm odd}}$ and $\psi=x^{-1}\phi\in\Cinf_{\text{\rm ev}}$, we get{\samepage
    \[
      |x|^\sigma|\phi(x)|\le
        \begin{cases}
          |\psi(x)| & \text{if $0<|x|\le1$},\\
          |\phi(x)| & \text{if $|x|\ge1$} .
        \end{cases}
    \]
  So $\|\phi\|_{\SS_{w,\sigma}^0}\le\max\{\|\phi\|_{\SS^0},\|\psi\|_{\SS^0}\}\le\|\phi\|_{\SS^1}$ by~\eqref{psi (m)(x)}.}

  Now, assume that $m\ge1$ and the result holds for $m-1$. Let $i,j\in\N$ such that $i+j\le m$, and let $\phi\in\SS_{\text{\rm ev}}\cup\SS_{\text{\rm odd}}$. Independently of the parity of $\phi$ and $i+j$, we have $|x|^\sigma |x^i\phi^{(j)}(x)|\le|x^i\phi^{(j)}(x)|$ for $|x|\ge1$.

  Suppose that $\phi\in\SS_{\text{\rm ev}}$. If $i=0$ and $j$ is odd, then $\phi^{(j)}\in\SS_{\text{\rm odd}}$. Thus there is some $\psi\in\SS_{\text{\rm ev}}$ such that $\phi^{(j)}=x\psi$, obtaining $|x|^\sigma|\phi^{(j)}(x)|\le|\psi(x)|$ for $0<|x|\le1$. If $i+j$ is odd and $i>0$, then $|x|^\sigma|x^i\phi^{(j)}(x)|\le|x^{i-1}\phi^{(j)}(x)|$ for $0<|x|\le1$. Hence, by~\eqref{psi (m)(x)}, there is some $C>0$, independent of $\phi$, such that
    \[
      \|\phi\|_{\SS_{w,\sigma}^m}\le C\max\{\|\phi\|_{\SS^m},\|\psi\|_{\SS^0}\}\le C\max\bigl\{\|\phi\|_{\SS^m},\|\phi^{(j)}\|_{\SS^1}\bigr\}\le C \|\phi\|_{\SS^{m+1}} .
    \]

  Finally, assume $\phi\in\SS_{\text{\rm odd}}$, and let $\psi=x^{-1}\phi\in\SS_{\text{\rm ev}}$. If $i$ is even and $j=0$, then $|x|^\sigma |x^i\phi(x)|\le|x^i\psi(x)|$ for $0<|x|\le1$. If $i+j$ is even and $j>0$, then
    \[
      |x|^\sigma\big|x^i\phi^{(j)}(x)\big|\le\big|x^i\psi^{(j)}(x)\big|+j |x|^\sigma\big|x^i\psi^{(j-1)}(x)\big|
    \]
  for $0<|x|\le1$ because $[\frac{{\rm d}^j}{{\rm d}x^j},x]=j \frac{{\rm d}^{j-1}}{{\rm d}x^{j-1}}$. Therefore, by~\eqref{psi (m)(x)} and the induction hypothesis, there are some $C',C''>0$, independent of $\phi$, such that
    \begin{gather*}
      \|\phi\|_{\SS_{w,\sigma}^m}\le C'\max\bigl\{\|\phi\|_{\SS^m},\|\psi\|_{\SS^m}+\|\psi\|_{\SS_{w,\sigma}^{m-1}}\bigr\}\le C'' \|\phi\|_{\SS^{m+1}} .
\tag*{\qed}
    \end{gather*}
\renewcommand{\qed}{}
\end{proof}

\begin{Lemma}\label{l:SS w,sigma m' subset SS m}
  If $\sigma<0$, then $\SS_{w,\sigma}^{m+1}\subset\SS^m$ continuously.
\end{Lemma}

\begin{proof}
  Let $i,j\in\N$ such that $i+j\le m$. Since
    \[
      \big|x^i\phi^{(j)}(x)\big|\le
        \begin{cases}
          |x|^\sigma\big|x^i\phi^{(j)}(x)\big| & \text{if $0<|x|\le1$},\\
          |x|^\sigma\big|x^{i+1}\phi^{(j)}(x)\big| & \text{if $|x|\ge1$} .
        \end{cases}
    \]
  for any $\phi\in\Cinf$, we get $\|\phi\|_{\SS^m}\le\|\phi\|_{\SS_{w,\sigma}^{m+1}}$.
\end{proof}

\begin{Corollary}\label{c:SS=SS w,sigma}
  $\SS=\SS_{w,\sigma}$ as Fr\'echet spaces.
\end{Corollary}

\begin{Corollary}\label{c:x^-1:SS_w,sigma,odd^m' to SS_w,sigma,ev^m}
  $x^{-1}$ defines a bounded operator $\SS_{w,\sigma,\text{\rm odd}}^{m'}\to\SS_{w,\sigma,\text{\rm ev}}^m$, where
    \[
      m'=
        \begin{cases}
          m+2+\frac{1}{2}\lceil\sigma\rceil(\lceil\sigma\rceil+3) & \text{if $\sigma\ge0$},\\
          m+3 & \text{if $\sigma<0$} .
        \end{cases}
    \]
\end{Corollary}

\begin{proof}
  If $\sigma\ge0$, the composite
    \[
      \begin{CD}
        \SS_{w,\sigma,\text{\rm odd}}^{m+2+\frac{1}{2}\lceil\sigma\rceil(\lceil\sigma\rceil+3)}\hookrightarrow\SS_{\text{\rm odd}}^{m+\lceil\sigma\rceil+1} @>{x^{-1}}>> \SS_{\text{\rm ev}}^{m+\lceil\sigma\rceil}\hookrightarrow\SS_{w,\sigma,\text{\rm ev}}^m
      \end{CD}
    \]
  is bounded by Lemmas~\ref{l:SS m+lceil sigma/2 rceil subset SS w,sigma m} and~\ref{l:SS_w,sigma^m' subset SS^m}. If $\sigma<0$, the composite
    \[
      \begin{CD}
        \SS_{w,\sigma,\text{\rm odd}}^{m+3}\hookrightarrow\SS_{\text{\rm odd}}^{m+2} @>{x^{-1}}>> \SS_{\text{\rm ev}}^{m+1}\hookrightarrow\SS_{w,\sigma,\text{\rm ev}}^m ,
      \end{CD}
    \]
  is bounded by Lemmas~\ref{l:SS m+1 subset SS w,sigma m} and~\ref{l:SS w,sigma m' subset SS m}.
\end{proof}

\begin{Corollary}\label{c:x^-1:SS_w,sigma,odd to SS_w,sigma,ev}
  $x^{-1}$ defines a continuous operator $\SS_{w,\sigma,\text{\rm odd}}\to\SS_{w,\sigma,\text{\rm ev}}$.
\end{Corollary}

\begin{Lemma}\label{l:SS_w,sigma^M_m subset SS_sigma^m}
  $\SS_{w,\sigma,\text{\rm ev/odd}}^{M_{m,\text{\rm ev/odd}}}\subset\SS_{\sigma,\text{\rm ev/odd}}^m$ continuously, where
    \begin{gather*}
      M_{m,\text{\rm ev/odd}} =
        \begin{cases}
          \dfrac{3m}{2}+\dfrac{m}{4} \lceil\sigma\rceil(\lceil\sigma\rceil+3) & \text{if $\sigma\ge0$ and $m$ is even},\\
          2m & \text{if $\sigma<0$ and $m$ is even} ,
        \end{cases}
      \\
      M_{m,\text{\rm ev}} =
        \begin{cases}
          \dfrac{3m-1}{2}+\dfrac{m-1}{4} \lceil\sigma\rceil(\lceil\sigma\rceil+3) & \text{if $\sigma\ge0$ and $m$ is odd},\\
          2m-1 & \text{if $\sigma<0$ and $m$ is odd} ,
        \end{cases}
      \\
      M_{m,\text{\rm odd}} =
        \begin{cases}
          \dfrac{3m+1}{2}+\dfrac{m+1}{4} \lceil\sigma\rceil(\lceil\sigma\rceil+3)
          & \text{if $\sigma\ge0$ and $m$ is odd},\\
          2m+1 & \text{if $\sigma<0$ and $m$ is odd} .
        \end{cases}
    \end{gather*}
\end{Lemma}

\begin{proof}
  This follows by induction on $m$. It is true for $m=0$ because $\SS_{w,\sigma}^0=\SS_\sigma^0$ as Banach spaces. Now, let $m\ge1$, and assume that the result holds for $m-1$.

  For $\phi\in\Cinf_{\text{\rm ev}}$, $i+j\le m$ with $j>0$, and $x\in\R$, we have $|x^i(T_\sigma^j\phi)(x)|=|x^i(T_\sigma^{j-1}\phi')(x)|$ with $\phi'\in \Cinf_{\text{\rm odd}}$, obtaining $\|\phi\|_{\SS_\sigma^m}\le\|\phi'\|_{\SS_\sigma^{m-1}}$. But, by the induction hypothesis and since $M_{m,\text{\rm ev}}=M_{m-1,\text{\rm odd}}+1$, there are some $C,C'>0$, independent of $\phi$, such that
    \[
      \|\phi'\|_{\SS_\sigma^{m-1}}\le C \|\phi'\|_{\SS_{w,\sigma}^{M_{m-1,\text{\rm odd}}}}\le C' \|\phi\|_{\SS_{w,\sigma}^{M_{m,\text{\rm ev}}}} .
    \]

  For $\phi\in \Cinf_{\text{\rm odd}}$, let $\psi=x^{-1}\phi$, and take $i$, $j$ and $x$ as above. Then
    \[
      \big|x^i(T_\sigma^j\phi)(x)\big|\le \big|x^i(T_\sigma^{j-1}\phi')(x)\big|
      +2 |\sigma|\big|x^i\big(T_\sigma^{j-1}\psi\big)(x)\big|
    \]
  with $\phi',\psi\in\Cinf_{\text{\rm ev}}$, obtaining $\|\phi\|_{\SS_\sigma^m}\le\|\phi'\|_{\SS_\sigma^{m-1}}+2 |\sigma| \|\psi\|_{\SS_\sigma^{m-1}}$. But, by the induction hypothesis, Corollary~\ref{c:x^-1:SS_w,sigma,odd^m' to SS_w,sigma,ev^m}, and since
    \[
      M_{m,\text{\rm odd}}=
        \begin{cases}
          M_{m-1,\text{\rm ev}}+2+\frac{1}{2}\lceil\sigma\rceil(\lceil\sigma\rceil+3) & \text{if $\sigma\ge0$},\\
          M_{m-1,\text{\rm ev}}+3 & \text{if $\sigma<0$} ,
        \end{cases}
    \]
  there are some $C,C'>0$, independent of $\phi$, such that
    \begin{gather*}
      \|\phi'\|_{\SS_\sigma^{m-1}}+2 |\sigma| \|\psi\|_{\SS_\sigma^{m-1}}
      \le C\Big(\|\phi'\|_{\SS_{w,\sigma}^{M_{m-1,\text{\rm ev}}}}
      +\|\psi\|_{\SS_{w,\sigma}^{M_{m-1,\text{\rm ev}}}}\Big)
      \le C' \|\phi\|_{\SS_{w,\sigma}^{M_{m,\text{\rm odd}}}} .\tag*{\qed}
    \end{gather*}
\renewcommand{\qed}{}
\end{proof}

\begin{Corollary}\label{c:SS_w,sigma subset SS_sigma}
  $\SS_{w,\sigma}\subset\SS_\sigma$ continuously.
\end{Corollary}

\begin{Corollary}\label{l:SS_ev/odd^M_m,ev/odd subset SS_sigma,ev/odd^m}
  $\SS_{\text{\rm ev/odd}}^{M'_{m,\text{\rm ev/odd}}}\subset\SS_{\sigma,\text{\rm ev/odd}}^m$ continuously, where, with the notation of Lemma~{\rm \ref{l:SS_w,sigma^M_m subset SS_sigma^m}},
    \[
      M'_{m,\text{\rm ev/odd}}=
        \begin{cases}
          M_{m,\text{\rm ev/odd}}+\lceil\sigma\rceil & \text{if $\sigma\ge0$},\\
          M_{m,\text{\rm ev/odd}}+1 & \text{if $\sigma<0$} .
        \end{cases}
    \]
\end{Corollary}

\begin{proof}
  This follows from Lemmas~\ref{l:SS m+lceil sigma/2 rceil subset SS w,sigma m},~\ref{l:SS m+1 subset SS w,sigma m} and~\ref{l:SS_w,sigma^M_m subset SS_sigma^m}.
\end{proof}

\section{Perturbed Sobolev spaces}\label{s:W sigma m}

Observe that $\SS_\sigma\subset L^2(\R,|x|^{2\sigma}{\rm d}x)$. Like in the case where $\SS$ is considered as domain, it is easy to check that, in $L^2(\R,|x|^{2\sigma}{\rm d}x)$, with domain $\SS_\sigma$, $B$ is adjoint of $B'$ and $L$ is symmetric.

\begin{Lemma}\label{l:SS_sigma is a core of LL}
   $\SS_\sigma$ is a core\footnote{Recall that a {\em core\/} of a closed densely def\/ined operator $T$ between Hilbert spaces is any subspace of its domain~$\DD(T)$ which is dense with the graph norm.} of~$\LL$.
\end{Lemma}

\begin{proof}
  Let $R$ denote the restriction of $\LL$ to $\SS_\sigma$. Then $\LL\subset\overline{R}\subset R^*\subset\LL^*=\LL$ in $L^2(\R,|x|^{2\sigma}{\rm d}x)$ because $\SS\subset\SS_\sigma$ by Corollaries~\ref{c:SS=SS w,sigma} and~\ref{c:SS_w,sigma subset SS_sigma}.
\end{proof}

For each $m\in\R$, let $W_\sigma^m$ be the Hilbert space completion of $\SS$ with respect to the scalar product $\langle\ ,\ \rangle_{W_\sigma^m}$ def\/ined by $\langle\phi,\psi\rangle_{W_\sigma^m}=\langle(1+\LL)^m\phi,\psi\rangle_\sigma$. The corresponding norm will be denoted by $\|\ \|_{W_\sigma^m}$,
whose equivalence class is independent of the parameter $s$ used to def\/ine~$L$. In particular,
$W_\sigma^0=L^2(\R,|x|^{2\sigma}{\rm d}x)$. As usual, $W_\sigma^{m'}\subset
W_\sigma^m$ when $m'>m$,  and let
$W_\sigma^\infty=\bigcap_mW_\sigma^m$ with the
induced Fr\'echet topology. Once more, there are direct sum decompositions into subspaces of even and odd (generalized) functions, $W_\sigma^m=W_{\sigma,\text{\rm ev}}^m\oplus W_{\sigma,\text{\rm odd}}^m$ ($m\in[0,\infty]$). By Lemma~\ref{l:SS_sigma is a core of LL}, $\SS_\sigma$ can be used instead of $\SS$ in the def\/inition of $W_\sigma^m$.

Obviously, $L$ def\/ines a bounded operator $W_\sigma^{m+2}\to W_\sigma^m$ for each $m$, and therefore a continuous operator on $W_\sigma^\infty$. By~\eqref{[L,Sigma]=B' Sigma+Sigma B'=0}, $\Sigma$ def\/ines a bounded operator on each $W_\sigma^m$, and therefore a~continuous operator on $W_\sigma^\infty$. Moreover $B$ and $B'$ def\/ine bounded operators $W_\sigma^{m+1}\to W_\sigma^m$ for each $m\in\N$; this follows easily by induction on $m$, using~\eqref{L} and~\eqref{[L,B]} (the details of the proof are omitted because this observation will not be used). Thus $L$, $\Sigma$, $B$ and $B'$ def\/ine bounded operators on~$W_\sigma^\infty$. Also on the spaces~$W_\sigma^m$, the parity of (generalized) functions is preserved by~$L$ and~$\Sigma$, and reversed by~$B$ and~$B'$. Observe that~$B'$ is not adjoint of~$B$ in~$W_\sigma^m$ for~$m\ne0$.

The motivation of our tour through perturbed Schwartz spaces is the following embedding results; the second one is a version of the Sobolev embedding theorem.

\begin{Proposition}\label{p:SS sigma,m' subset W sigma m}
  $\SS_\sigma^{m+1}\subset W_\sigma^m$ continuously for all $m\in\N$.
\end{Proposition}

\begin{Proposition}\label{p:Sobolev}
  For all $m\in\N$, $W_\sigma^{m'}\subset\SS_\sigma^m$ continuously if $m'-m>1$.
\end{Proposition}

\begin{Corollary}\label{c:Sobolev}
  $\SS_\sigma=W_\sigma^\infty$ as Fr\'echet spaces.
\end{Corollary}

For each non-commutative polynomial $p$ (of two variables, $X$ and $Y$), let $p'$ denote the non-commutative polynomial obtained by reversing the order of the variables in $p$; e.g., if $p(X,Y)=X^2Y^3X$, then $p'(X,Y)=XY^3X^2$. It will be said that $p$ is {\em symmetric\/} if $p(X,Y)=p'(Y,X)$. Note that $p'(Y,X)p(X,Y)$ is symmetric for any $p$. Given any non-commutative polynomial $p$, the continuous operators $p(B,B')$ and $p'(B',B)$ on $\SS_\sigma$ are adjoint of each other in $L^2(\R,|x|^{2\sigma}{\rm d}x)$; thus $p(B,B')$ is a symmetric operator if $p$ is symmetric.

\begin{Lemma}\label{l:(1+L) m}
  For $m\in\N$, we have $(1+L)^m=\sum_aq'_a(B',B)q_a(B,B')$ for some finite family of homogeneous non-commutative polynomials $q_a$ of degree $\le m$.
\end{Lemma}

\begin{proof}
  The result follows easily from the following assertions.

    \begin{claim}\label{cl:L m, m even}
      If $m$ is even, then $L^m=g_m(B,B')^2$ for some symmetric homogeneous non-commuta\-ti\-ve polynomial $g_m$ of degree $m$.
    \end{claim}

    \begin{claim}\label{cl:L m, m odd}
      If $m$ is odd, then
        \[
          L^m=g_{m,1}'(B',B)g_{m,1}(B,B')+g_{m,2}'(B',B)g_{m,2}(B,B')
        \]
      for some homogeneous non-commutative polynomials $g_{m,1}$ and $g_{m,2}$ of degree $m$.
    \end{claim}

  If $m$ is even, then $L^{m/2}=g_m(B,B')$ for some
symmetric homogeneous non-commutative polynomial $g_m$ of degree $m$
by~\eqref{L}. So $L^m=g_m(B,B')^2$, showing Claim~\ref{cl:L m, m even}.

  If $m$ is odd, write $L^{\lfloor m/2\rfloor}=f_m(B,B')$ as above for some symmetric homogeneous non-commutative polynomial $f_m$ of degree $m-1$. Then, by~\eqref{L},
    \[
      L^m=\frac{1}{2} f_m(B,B')(BB'+B'B)f_m(B,B') .
    \]
  Thus Claim~\ref{cl:L m, m odd} follows with
    \begin{gather*}
      g_{m,1}(B,B')=\frac{1}{\sqrt{2}} B'f_m(B,B') ,\qquad
      g_{m,2}(B,B')=\frac{1}{\sqrt{2}} Bf_m(B,B') .\tag*{\qed}
    \end{gather*}
\renewcommand{\qed}{}
\end{proof}

\begin{proof}[Proof of Proposition~\ref{p:SS sigma,m' subset W sigma m}]
  By the def\/initions of $B$ and $B'$, and~\eqref{SS_sigma^m, sigma ge0}--\eqref{SS_sigma,odd^m, sigma<0}, for each homogeneous non-commutative polynomial $p$ of three variables with degree $d\le m+1$, there exists some $C_p>0$ such that, for all $\phi\in\SS_\sigma$: if $\sigma<0$, $|x|\le1$, and $\phi$ and $d$ have the same parity, then $|(p(x,B,B')\phi)(x)|\le C_p \|\phi\|_{\SS_\sigma^{m+1}}$; and, otherwise, $|x|^\sigma |(p(x,B,B')\phi)(x)|\le C_p \|\phi\|_{\SS_\sigma^{m+1}}$.

  With the notation of Lemma~\ref{l:(1+L) m}, let $d_a$ denote the degree of each $q_a$, and let $\bar q_a(x,B,B')=x q_a(B,B')$. If $\sigma\ge0$, then
    \begin{gather*}
      \|\phi\|_{W_\sigma^m}^2=\sum_a\|q_a(B,B')\phi\|_\sigma^2
      =\sum_a\int_{-\infty}^\infty|(q_a(B,B')\phi)(x)|^2 |x|^{2\sigma}{\rm d}x\\
\hphantom{\|\phi\|_{W_\sigma^m}^2=}{} \le 2\sum_a\left(C_{q_a}^2+C_{\bar q_a}^2\int_1^\infty x^{-2}{\rm d}x\right)\|\phi\|_{\SS_\sigma^{m+1}}^2
    \end{gather*}
for $\phi\in\SS_\sigma$. Similarly, if $\sigma<0$, then $\|\phi\|_{W_\sigma^m}^2$ is bounded by
    \[
      2\left(\sum_{\text{\rm$d_a$ even}}C_{q_a}^2\int_0^1x^{2\sigma}{\rm d}x
      +\sum_{\text{\rm$d_a$ odd}}C_{q_a}^2+\sum_aC_{\bar q_a}^2\int_1^\infty x^{-2}{\rm d}x\right)\|\phi\|_{\SS_\sigma^{m+1}}^2
    \]
   for $\phi\in\SS_{\sigma,\text{\rm ev}}$, and by
     \[
      2\left(\sum_{\text{\rm$d_a$ even}}C_{q_a}^2+\sum_{\text{\rm$d_a$ odd}}C_{q_a}^2\int_0^1x^{2\sigma}{\rm d}x
      +\sum_aC_{\bar q_a}^2\int_1^\infty x^{-2}{\rm d}x\right)\|\phi\|_{\SS_\sigma^{m+1}}^2
    \]
  for $\phi\in\SS_{\sigma,\text{\rm odd}}$.
\end{proof}

For $c=(c_k)$ and $c'=(c'_k)$ in $\R^{\N}$, and $m\in\R$, the expressions $\|c\|_{\CC_m}=\sup_k|c_k|(1+k)^m$ and $\langle c,c'\rangle_{\ell^2_m}=\sum_kc_kc'_k(1+k)^m$ def\/ine a norm and scalar product with possible inf\/inite values, and let $\|\ \|_{\ell^2_m}$ denote the norm with possible inf\/inite values induced by $\langle\ ,\ \rangle_{\ell^2_m}$. Then $\CC_m=\{c\in\R^{\N}\mid\|c\|_{\CC_m}<\infty\}$ becomes a Banach space with $\|\ \|_{\CC_m}$, and $\ell^2_m=\{c\in\R^{\N}\mid\|c\|_{\ell^2_m}<\infty\}$ is a Hilbert space with $\langle\ ,\ \rangle_{\ell^2_m}$. If $m'>m$, then $\CC_{m'}\subset\CC_m$ and $\ell^2_{m'}\subset\ell^2_m$ continuously. Let $\CC_\infty=\bigcap_m\CC_m$ and $\ell^2_\infty=\bigcap_m\ell^2_m$, with the induced Fr\'echet topologies. It is said that $c$ is even/odd if $c_k=0$ for odd/even $k$. There are direct sum decompositions into subspaces of even and odd sequences, $\CC_m=\CC_{m,\text{\rm ev}}\oplus\CC_{m,\text{\rm odd}}$ and $\ell^2_m=\ell^2_{m,\text{\rm ev}}\oplus\ell^2_{m,\text{\rm odd}}$ ($m\in\R\cup\{\infty\}$).

\begin{Lemma}\label{l:ell 2 2m subset CC m}
  $\ell^2_{2m}\subset\CC_m$ and $\CC_{m'}\subset\ell^2_m$ continuously if $2m'-m>1$.
\end{Lemma}

\begin{proof}
  It is easy to see that $\|c\|_{\CC_m}\le\|c\|_{\ell^2_{2m}}$ and $ \|c\|_{\ell^2_m}\le\|c\|_{\CC_{m'}}\big(\sum_k(1+k)^{m-2m'}\big)^{1/2}$ for any $c\in\CC_\infty$, where the last series is convergent because $m-2m'<-1$.
\end{proof}

\begin{Corollary}\label{c:ell 2 infty=CC}
  $\ell^2_\infty=\CC_\infty$ as Fr\'echet spaces.
\end{Corollary}

According to Section~\ref{ss:L}, the ``Fourier coef\/f\/icients'' mapping $\phi\mapsto(\langle\phi_k,\phi\rangle_\sigma)$ def\/ines a quasi-isometry $W_\sigma^m\to\ell^2_m$ for all f\/inite $m$, and therefore an isomorphism $W_\sigma^\infty\to\CC_\infty$ of Fr\'echet espaces. This map is compatible with the decompositions into even and odd subspaces.

\begin{Corollary}\label{c:SS_sigma, rapidly degreasing Fourier coefficients}
  Any $\phi\in L^2\big(\R,|x|^{2\sigma}{\rm d}x\big)$ is in $\SS_\sigma$ if and only if its ``Fourier coefficients'' $\langle\phi_k,\phi\rangle_\sigma$ are rapidly decreasing on $k$.
\end{Corollary}

\begin{proof}
  By Corollary~\ref{c:Sobolev}, the ``Fourier coef\/f\/icients'' mapping def\/ines an isomorphism $\SS_\sigma\to\CC_\infty$ of Fr\'echet spaces.
\end{proof}

The operator $\ell^2_{m'}\hookrightarrow\ell^2_m$ is compact if $m'>m$ (see e.g.\ \cite[Theorem 5.8]{Roe1998}). So, by using  the ``Fourier coef\/f\/icients'' mapping, the operator $W_\sigma^{m'}\hookrightarrow W_\sigma^m$ is compact if $m'>m$ (a version of the Rellich theorem).

\begin{proof}[Proof of Proposition~\ref{p:Sobolev}]
  For $\phi\in\SS_\sigma$, its ``Fourier coef\/f\/icients'' $c_k=\langle\phi_k,\phi\rangle_\sigma$ form a sequence $c=(c_k)$ in $\CC_\infty$, and $\sum_k|c_k| (1+k)^{m/2}\le\|c\|_{\ell^2_{m'}}(\sum_k(1+k)^{m-m'})^{1/2}$ by Cauchy--Schwartz inequality, where the last series is convergent since $m-m'<-1$. Therefore there is some $C>0$, independent of $\phi$, such that
    \begin{gather}\label{sum k|c k|(1+k) m/2}
      \sum_k|c_k| (1+k)^{m/2}\le C \|\phi\|_{W_\sigma^{m'}} .
    \end{gather}

  On the other hand, for all $i,j\in\N$ with $i+j\le m$, there is some homogeneous non-commutative polynomial $p_{ij}$ of degree $i+j$ such that $x^iT_\sigma^j=p_{ij}(B,B')$. Then, by~\eqref{phi k}--\eqref{B phi k}, there is some $C_{ij}>0$, independent of $\phi$, such that
    \begin{gather}\label{|langle phi k,x iD sigma j phi rangle sigma|}
      |\langle\phi_k,x^iT_\sigma^j\phi\rangle_\sigma|\le C_{ij}(1+k)^{m/2}\sum_{|\ell-k|\le m}|c_\ell| .
    \end{gather}

 Assume that $\sigma\ge0$. By~\eqref{sum k|c k|(1+k) m/2},~\eqref{|langle phi k,x iD sigma j phi rangle sigma|} and Theorem~\ref{t:upper estimates of xi k}(i), there is some $C'_{ij}>0$, independent of $\phi$ and $x_0$, so that
    \begin{gather}
      |x_0|^\sigma |(x^iT_\sigma^j\phi)(x_0)|\le |x_0|^\sigma\sum_k|\langle\phi_k,x^iT_\sigma^j\phi\rangle_\sigma| |\phi_k(x_0)|\nonumber\\
 \hphantom{|x_0|^\sigma |(x^iT_\sigma^j\phi)(x_0)|}{}
      =\sum_k|\langle\phi_k,x^iT_\sigma^j\phi\rangle_\sigma| |\xi_k(x_0)|\le C'_{ij} \|\phi\|_{W_\sigma^{m'}}\label{|x| sigma/2 |x iD sigma j phi(x)|}
    \end{gather}
  for all $x_0\in\R$. Hence $\|\phi\|_{\SS_\sigma^m}\le C'\|\phi\|_{W_\sigma^{m'}}$ for some $C'>0$ independent of $\phi$.

  Now, suppose that  $\sigma<0$. From~\eqref{sum k|c k|(1+k) m/2},~\eqref{|langle phi k,x iD sigma j phi rangle sigma|}, and Theorem~\ref{t:upper estimates of xi k}(ii),(iii), it follows that there is some $C'_{ij}>0$, independent of $\phi$ and $x_0$, such that
    \[
      |(x^iT_\sigma^j\phi)(x_0)|\le\sum_k|\langle\phi_k,x^iT_\sigma^j\phi\rangle_\sigma| |\phi_k(x_0)|
      \le C'_{ij} \|\phi\|_{W_\sigma^{m'}}
    \]
  if $\phi$ and $i+j$ have the same parity, and $|x_0|\le1$; otherwise, $|x_0|^\sigma |(x^iT_\sigma^j\phi)(x_0)|\le C'_{ij} \|\phi\|_{W_\sigma^{m'}}$ like in~\eqref{|x| sigma/2 |x iD sigma j phi(x)|}. So $\|\phi\|_{\SS_ \sigma ^m}\le C'\|\phi\|_{W_\sigma^{m'}}$ with $C'>0$ independent of $\phi$.
\end{proof}

As suggested by~\eqref{x -1 p k}, consider the mapping $c=(c_k)\mapsto\Xi(c)=(d_\ell)$, where $c$ is odd and $\Xi(c)$ is even, with
  \[
    d_\ell=\sum_{k\in\{\ell+1,\ell+3,\dots\}}(-1)^{\frac{k-\ell-1}{2}}\sqrt{\frac{(k-1)(k-3)\cdots(\ell+2)2s}{(k+2\sigma)(k-2+2\sigma)\cdots(\ell+1+2\sigma)}}
    c_k
  \]
for $\ell$ even, assuming that this series is convergent.

\begin{Lemma}\label{l:Xi}
  $\Xi$ defines a bounded map $\ell^2_{m',\text{\rm odd}}\to\CC_{m,\text{\rm ev}}$ if $m'-m>1$.
\end{Lemma}

\begin{proof}
  By the Cauchy--Schwartz  inequality,
    \begin{gather*}
      \|d\|_{\CC_m} =\sup_\ell\sum_{k\in\{\ell+1,\ell+3,\dots\}}\sqrt{\frac{(k-1)(k-3)\cdots(\ell+2)2s}{(k+2\sigma)(k-2+2\sigma)\cdots(\ell+1+2\sigma)}} |c_k| (1+\ell)^m\\
\hphantom{\|d\|_{\CC_m}}{}
\le\sqrt{2s}\sup_\ell\sum_{k\in\{\ell+1,\ell+3,\dots\}}|c_k| (1+\ell)^m\\
\hphantom{\|d\|_{\CC_m}}{}
\le\sqrt{2s} \|c\|_{\ell^2_{m'}}\sup_\ell\biggl(\sum_{k\in\{\ell+1,\ell+3,\dots\}}(1+k)^{-m'}(1+\ell)^m\biggr)^{1/2}\\
\hphantom{\|d\|_{\CC_m}}{}
\le\sqrt{2s} \|c\|_{\ell^2_{m'}}\biggl(\sum_k(1+k)^{m-m'}\biggr)^{1/2} ,
    \end{gather*}
  where the last series is convergent because $m-m'<-1$.
\end{proof}

 \begin{Corollary}\label{c:x^-1:SS_sigma,odd^m_1 to SS_sigma,ev^m bounded}
  $x^{-1}$ defines a bounded operator $\SS_{\sigma,\text{\rm odd}}^{m'}\to\SS_{\sigma,\text{\rm ev}}^m$ if $2m'>m+6$.
\end{Corollary}

\begin{proof}
  Since $2m'>m+6$, there are $m_1,m_2\in\R$ so that $m'-m_2>2$, $2m_2-m_1>1$ and $m_1-m>1$. Then, by Propositions~\ref{p:SS sigma,m' subset W sigma m} and~\ref{p:Sobolev}, Lemmas~\ref{l:ell 2 2m subset CC m} and~\ref{l:Xi}, and using the ``Fourier coef\/f\/icients'' mapping, we get the composition of bounded maps,
    \[
      \SS_{\sigma,\text{\rm odd}}^{m'}\hookrightarrow W_{\sigma,\text{\rm odd}}^{m'-1}\rightarrow\ell^2_{m'-1,\text{\rm odd}}\overset{\Xi}{\longrightarrow}\CC_{m_2,\text{\rm ev}}\hookrightarrow\ell^2_{m_1,\text{\rm ev}}\rightarrow W_{\sigma,\text{\rm ev}}^{m_1}\hookrightarrow\SS_{\sigma,\text{\rm ev}}^m ,
    \]
  which extends $x^{-1}:\SS_{\text{\rm odd}}\to\SS_{\text{\rm ev}}$ by~\eqref{x -1 p k}.
\end{proof}

\begin{quest}
  Is it possible to prove Corollary~\ref{c:x^-1:SS_sigma,odd^m_1 to SS_sigma,ev^m bounded} without using~\eqref{x -1 p k}?
\end{quest}

\begin{Corollary}\label{c:x^-1 continuous}
  $x^{-1}$ defines a continuous operator $\SS_{\sigma,\text{\rm odd}}\to\SS_{\sigma,\text{\rm ev}}$.
\end{Corollary}

\begin{Lemma}\label{l:SS_sigma^m+2 subset SS_w,sigma^m}
  $\SS_{\sigma,\text{\rm ev/odd}}^{M_{m,\text{\rm ev/odd}}}\subset\SS_{w,\sigma,\text{\rm ev/odd}}^m$ continuously for all $m$, where $M_{0,\text{\rm ev/odd}}=0$, $M_{1,\text{\rm ev}}=1$, $M_{1,\text{\rm odd}}=M_{2,\text{\rm odd}}=4$, $M_{2,\text{\rm ev}}=M_{3,\text{\rm ev}}=5$, $M_{3,\text{\rm odd}}=6$, and $M_{m,\text{\rm ev/odd}}=m+2$ for $m\ge4$.
\end{Lemma}

\begin{proof}
  We proceed by induction on $m$. The case $m=0$ was already indicated in the proof of Lemma~\ref{l:SS_w,sigma^M_m subset SS_sigma^m}. Now, let $m\ge1$ and assume that the result holds for $m-1$.

  For $\phi\in \Cinf_{\text{\rm ev}}$, $i+j\le m$ with $j>0$ and $x\in\R$, we have $|x^i\phi^{(j)}(x)|=|x^i(T_\sigma\phi)^{(j-1)}(x)|$ with $T_\sigma\phi\in \Cinf_{\text{\rm odd}}$, obtaining $\|\phi\|_{\SS_{w,\sigma}^m}\le\|T_\sigma\phi\|_{\SS_{w,\sigma}^{m-1}}$. But, by the induction hypothesis and because $M_{m,\text{\rm ev}}=M_{m-1,\text{\rm odd}}+1$, there are some $C,C'>0$, independent of $\phi$, such that $\|T_\sigma\phi\|_{\SS_{w,\sigma}^{m-1}}\le C \|T_\sigma\phi\|_{\SS_\sigma^{M_{m-1,\text{\rm odd}}}}\le C' \|\phi\|_{\SS_\sigma^{M_{m,\text{\rm ev}}}}$.

  For $\phi\in \Cinf_{\text{\rm odd}}$, let $\psi=x^{-1}\phi$, and take $i$, $j$ and $x$ as above. We have
    \[
      \big|x^i\phi^{(j)}(x)\big|\le\big|x^i(T_\sigma\phi)^{(j-1)}(x)\big|
      +2 |\sigma|\big|x^i\psi^{(j-1)}(x)\big|
    \]
  with $T_\sigma\phi,\psi\in\Cinf_{\text{\rm ev}}$, obtaining $ \|\phi\|_{\SS_{w,\sigma}^m}\le\|T_\sigma\phi\|_{\SS_{w,\sigma}^{m-1}}
      +2 |\sigma| \|\psi\|_{\SS_{w,\sigma}^{m-1}}$. But, by the induction hypothesis, Corollary~\ref{c:x^-1:SS_sigma,odd^m_1 to SS_sigma,ev^m bounded}, and since $M_{m,\text{\rm odd}}\ge M_{m-1,\text{\rm ev}}+1$ and $2M_{m,\text{\rm odd}}>M_{m-1,\text{\rm ev}}+6$, there are some $C,C'>0$, independent of $\phi$, such that
    \begin{gather*}
      \|T_\sigma\phi\|_{\SS_{w,\sigma}^{m-1}}+2 |\sigma| \|\psi\|_{\SS_{w,\sigma}^{m-1}}
      \le C\left(\|\phi'\|_{\SS_\sigma^{M_{m-1,\text{\rm ev}}}}+\|\psi\|_{\SS_\sigma^{M_{m-1,\text{\rm ev}}}}\right)
      \le C' \|\phi\|_{\SS_\sigma^{M_{m,\text{\rm odd}}}} .\tag*{\qed}
    \end{gather*}
\renewcommand{\qed}{}
\end{proof}

\begin{Corollary}\label{c:SS_sigma^m' subset SS^m}
  $\SS_{\sigma,\text{\rm ev/odd}}^{M'_{m,\text{\rm ev/odd}}}\subset\SS_{\text{\rm ev/odd}}^m$ continuously, where, with the notation of Lemma~{\rm \ref{l:SS_sigma^m+2 subset SS_w,sigma^m}}, $M'_{m,\text{\rm ev/odd}}=M_{m_\sigma,\text{\rm ev/odd}}$ for $m_\sigma=m+1+\frac{1}{2}\lceil\sigma\rceil(\lceil\sigma\rceil+1)$.
\end{Corollary}

\begin{proof}
  This follows from Lemmas~\ref{l:SS_w,sigma^m' subset SS^m},~\ref{l:SS w,sigma m' subset SS m} and~\ref{l:SS_sigma^m+2 subset SS_w,sigma^m}.
\end{proof}

\begin{Corollary}\label{c:SS_sigma = SS}
  $\SS_\sigma=\SS$ as Fr\'echet spaces.
\end{Corollary}

\begin{proof}
  This is a consequence of Corollaries~\ref{l:SS_ev/odd^M_m,ev/odd subset SS_sigma,ev/odd^m} and~\ref{c:SS_sigma^m' subset SS^m}
\end{proof}

Corollaries~\ref{l:SS_ev/odd^M_m,ev/odd subset SS_sigma,ev/odd^m} and~\ref{c:SS_sigma^m' subset SS^m} and Propositions~\ref{p:SS sigma,m' subset W sigma m} and~\ref{p:Sobolev} give Theorems~\ref{t:SS^m' subset W_sigma^m} and~\ref{t:Sobolev}.

\section[Perturbations of $H$ on $\R_+$]{Perturbations of $\boldsymbol{H}$ on $\boldsymbol{\R_+}$}\label{s:P}

Since the function $|x|^{2\sigma}$ is even, the decomposition $\SS=\SS_{\text{\rm ev}}\oplus\SS_{\text{\rm odd}}$ extends to an orthogonal decomposition
\[
    L^2\big(\R,|x|^{2\sigma}{\rm d}x\big)
    =L^2_{\text{\rm ev}}\big(\R,|x|^{2\sigma}{\rm d}x\big)\oplus L^2_{\text{\rm odd}}\big(\R,|x|^{2\sigma}{\rm d}x\big) .
\]
Let $L_{\text{\rm ev/odd}}$ and $\LL_{\text{\rm ev/odd}}$, or $L_{\sigma,\text{\rm ev/odd}}$ and $\LL_{\sigma,\text{\rm ev/odd}}$, denote the corresponding components of~$L$ and~$\LL$. $L_{\text{\rm ev/odd}}$ is essentially self-adjoint in $L^2_{\text{\rm ev/odd}}(\R,|x|^{2\sigma}{\rm d}x)$, and its self-adjoint extension is~$\LL_{\text{\rm ev/odd}}$, which satisf\/ies an obvious version of Corollary~\ref{c:SS=W_sigma^infty}.

Fix an open subset $U\subset\R_+$ of full Lebesgue measure. Let $\SS_{\text{\rm ev/odd},U}\subset\Cinf(U)$ be the linear subspace of restrictions to $U$ of the functions in $\SS_{\text{\rm ev/odd}}$. The restriction to $U$ def\/ines a linear isomorphism $\SS_{\text{\rm ev/odd}}\cong\SS_{\text{\rm ev/odd},U}$, and a unitary isomorphism $L^2_{\text{\rm ev/odd}}(\R,|x|^{2\sigma}{\rm d}x)\cong L^2(U,2x^{2\sigma}{\rm d}x)$. Via these isomorphisms, $L_{\text{\rm ev/odd}}$ corresponds to an operator $L_{\text{\rm ev/odd},U}$ on $\SS_{\text{\rm ev/odd},U}$, and $\LL_{\text{\rm ev/odd}}$ corresponds to a self-adjoint operator $\LL_{\text{\rm ev/odd},U}$ in $L^2(U,x^{2\sigma}{\rm d}x)$; the more explicit notation $L_{\sigma,\text{\rm ev/odd},U}$ and $\LL_{\sigma,\text{\rm ev/odd},U}$ may be used. Let $\phi_{k,U}=\phi_k|_U$, whose norm in $L^2(U,x^{2\sigma}{\rm d}x)$ is $1/\sqrt{2}$.

Going one step further, for any positive function $h\in C^2(U)$, the multiplication by $h$ def\/ines a unitary isomorphism $h:L^2(U,x^{2\sigma}{\rm d}x)\to L^2(U,x^{2\sigma}h^{-2}{\rm d}x)$. Thus $hL_{\text{\rm ev},U}h^{-1}$, with domain $h \SS_{\text{\rm ev},U}$, is essentially self-adjoint in $L^2(U,x^{2\sigma}h^{-2}{\rm d}x)$, and its self-adjoint extension is $h\LL_{\text{\rm ev},+}h^{-1}$. Via these unitary isomorphisms, we get an obvious version of Corollary~\ref{c:SS=W_sigma^infty} for $h\LL_{\text{\rm ev},+}h^{-1}$. By using
  \begin{gather}\label{[frac partial partial x,h]}
    \left[\frac{{\rm d}}{{\rm d}x},h\right]=h' ,\qquad
    \left[\frac{{\rm d}^2}{{\rm d}x^2},h\right]=2h'\frac{{\rm d}}{{\rm d}x}+h'' ,
  \end{gather}
it easily follows that $hL_{\text{\rm ev},U}h^{-1}$has the form of $P$ in Theorem~\ref{t:P}. Then Theorem~\ref{t:P} is a~consequence of the following.

\begin{Lemma}\label{l:P}
  For $\sigma>-1/2$, a positive function $h\in C^2(U)$, and $P=H-2f_1 \frac{{\rm d}}{{\rm d}x}+f_2$ with $f_1\in C^1(U)$ and $f_2\in C(U)$, we have $P=hL_{\sigma,\text{\rm ev},U}h^{-1}$ on $h \SS_{\text{\rm ev},U}$ if and only if~$f_1$, $f_2$ and $h$ satisfy the conditions of Theorem~{\rm \ref{t:P}}.
\end{Lemma}

\begin{proof}
  By~\eqref{[frac partial partial x,h]},
    \[
      h^{-1}Ph=H-2\big(h^{-1}h'+f_1\big)\frac{{\rm d}}{{\rm d}x}-h^{-1}h''-2h^{-1}f_1h'+f_2 .
    \]
  So $P=hL_{\sigma,\text{\rm ev},U}h^{-1}$ if and only if $h^{-1}h'=\sigma x^{-1}-f_1$ and $f_2=h^{-1}h''+2h^{-1}h'f_1$, which are easily seen to be equivalent to the conditions of Theorem~\ref{t:P}.
\end{proof}

\begin{Remark}\label{r:Q}
  By~\eqref{[frac partial partial x,h]}, we get an operator of the same type if $h$ and $\frac{{\rm d}}{{\rm d}x}$ is interchanged in the operator $P$ of Theorem~\ref{t:P}.
\end{Remark}

\begin{Remark}
By using~\eqref{[frac partial partial x,h]} with $h=x^{-1}$ on $\R_+$, it is easy to check that $L_{\sigma,\text{\rm odd},\R_+}=xL_{1+\sigma,\text{\rm ev},\R_+}x^{-1}$ on $\SS_{\text{\rm odd},\R_+}=x \SS_{\text{\rm ev},\R_+}$. So no new operators are obtained with this process by using $L_{\sigma,\text{\rm odd}}$ instead of $L_{\sigma,\text{\rm ev}}$.
\end{Remark}

\subsection*{Acknowledgements}

The f\/irst author is partially supported by MICINN and MEC, Grants MTM2008-02640,\linebreak  MTM2011-25656  and PR2009-0409.

\pdfbookmark[1]{References}{ref}
\LastPageEnding

\end{document}